\documentclass[a4paper]{article}
\usepackage[utf8x]{inputenc}
\usepackage{t1enc}
\PrerenderUnicode{ĂăÁáÄäÂâČčĆćĐđÉéĚěËëÍíÎîÏïĹĺÔôÖöÓóŐőßŠšÛûŽž}
\usepackage{latexsym}
\usepackage{amsfonts}
\usepackage{amssymb}
\usepackage{amsbsy}
\usepackage{amsmath}
\usepackage{paralist}
\usepackage{graphicx}
\usepackage{subfig}
\usepackage{booktabs}
\usepackage{multicol}
\usepackage[ruled,lined,algosection,longend]{algorithm2e}
\usepackage{arydshln}
\usepackage{ifthen}

\usepackage{url}

\numberwithin{equation}{section}

\newcommand{\diag}{\operatorname{diag}}
\newcommand{\assgn}{\ensuremath\mathrel{\mathop:}=}
\newcommand{\bor}{\mathrel{\operatorname{bit\_or}}}
\newcommand\hpm{\hphantom{-}}
\newcommand\hpz{\hphantom{0}}
\newcommand{\id}{\operatorname{id}}
\newcommand{\fl}{\operatorname{fl}}
\newcommand{\hyp}{{\it hyp}}
\newcommand{\Teps}{T_{\varepsilon}}
\newcommand{\llrightarrow}{\relbar\joinrel\relbar\joinrel\relbar\joinrel\relbar\joinrel\rightarrow}
\newcommand\iblk{\textit{i\_blk}}
\newcommand\jblk{\textit{j\_blk}}
\newcommand\ip{\textit{ip}}
\newcommand\jp{\textit{jp}}

\begin{document}
\SetAlFnt{\small}

\date{\today}

\title{A GPU-based hyperbolic SVD algorithm%
  \thanks{This work was supported by grants 037--1193086--2771 and
    120--1201703--1672 by the Ministry of Science, Education and Sports,
    Republic of Croatia.%
  }%
}
\author{%
Vedran Novaković%
  \thanks{Faculty of Mechanical Engineering and Naval Architecture,
  University of Zagreb, Ivana Lučića 5, 10000 Zagreb, Croatia,
  e-mail: {\tt venovako@fsb.hr}},\hbox{\ }
Sanja Singer%
  \thanks{Faculty of Mechanical Engineering and Naval Architecture,
  University of Zagreb, Ivana Lučića 5, 10000, Croatia,
  e-mail: {\tt ssinger@fsb.hr}}}

\maketitle

\markboth{V.\ Novaković and S.\ Singer}
  {A GPU-based hyperbolic SVD algorithm}

\begin{abstract}
A one-sided Jacobi hyperbolic singular value decomposition (HSVD)
algorithm, using a massively parallel graphics processing unit (GPU),
is developed.
The algorithm also serves as the final stage of solving a symmetric
indefinite eigenvalue problem.  Numerical testing demonstrates the
gains in speed and accuracy over sequential and MPI-parallelized
variants of similar Jacobi-type HSVD algorithms.  Finally,
possibilities of hybrid CPU--GPU parallelism are discussed.
\end{abstract}

%%%%% Keywords %%%%%%%%%%%

\vspace*{6pt}

%%%%% AMS classification %%%%%%%%%%%

\vspace*{6pt}
%
%%%%%%%%%%%%%%%%%%%%%%%%%%%%%%%%%%%%%%%%%%%%%%%%%%%%%%%%%%%%%%%%%%%%%%%%%%%%%%
%
\section{Introduction}
\label{sec:1}
%
%%%%%%%%%%%%%%%%%%%%%%%%%%%%%%%%%%%%%%%%%%%%%%%%%%%%%%%%%%%%%%%%%%%%%%%%%%%%%%
%
In this paper a hyperbolic SVD algorithm (HSVD) for graphics
processors (GPUs) is developed.  To the best of our knowledge, this is 
the first one-sided Jacobi-type HSVD algorithm for full column rank
matrices of arbitrary dimensions, using GPUs.

The first paper that we know of, which describes a Jacobi-type
computation of SVD on GPUs, was published by Zhang and Dou in
\cite{Zhang-Dou-07}.  According to the English summary of the paper,
they compute the SVD by using the one-sided Jacobi algorithm on a
matrix of order $512$.  The order of orthogonalization of pairs of
columns is chosen by the parallel pivot strategy described in
\cite{Brent-LukF-85}.  Later, Lahabar and Narayanan
\cite{Lahabar-Narayanan-09} have computed SVD on GPUs by
bidiagonalization followed by a bidiagonal QR algorithm (in single
precision arithmetic).  Finally, Sachdev, Vanjani and Hall in
\cite{Sachdev-Vanjani-Hall-10} have presented an algorithm for Takagi
SVD (for complex symmetric matrices) by a symmetrized version of the
two-sided Kogbetliantz algorithm (in single precision arithmetic).

Given a Hermitian indefinite matrix $M$, factorized as
$M = G J G^{\ast}$, with $G$ of full column rank and $J$ a diagonal
matrix holding the inertia of $M$, we aim at computing the HSVD of the
factor $G$,
\begin{equation}
  G = U \begin{bmatrix}
    \Sigma \\
    0
  \end{bmatrix} V^{\ast},
\label{1.1}
\end{equation}
where $U$ is unitary, $V$ is $J$-unitary (i.e., $V^{\ast} J V = J$),
and $\Sigma$ is a diagonal matrix with positive diagonal entries.

If $G$ is not of the full column rank, $\Sigma$ in (\ref{1.1}),
is not diagonal \cite[Remark 5]{Zha-96}, and the non-diagonal part
reflects the loss of rank of $G J G^{\ast}$ compared with the rank of
$G$.  In the (very important) special case, when $J = I$, the HSVD is
the ordinary SVD.

If the HSVD of a factor $G$ is known, then the eigendecomposition of
$M$ is readily available as
\begin{displaymath}
  M = G J G^{\ast}
    = U \begin{bmatrix}\Sigma\\0\end{bmatrix} V^{\ast} J
      V\begin{bmatrix}\Sigma & 0\end{bmatrix}U^{\ast}
    = U (\Sigma^2 J) U^{\ast},
\end{displaymath}
i.e., the matrix $U$ of left singular vectors of $G$ is also the
eigenvector matrix of $M$, while the eigenvalues of $M$ are diagonal
elements of $\Sigma^2 J$.

In many applications, $M$ is already given implicitly, by its
rectangular factor $G$, and the signature matrix $J$.  That being the
case, $G$ should be shortened either by the QR factorization (if $G$
is tall) or by the hyperbolic QR factorization (see
\cite{SingerSanja-06}) of $G^{\ast}$ (if $G^{\ast}$ is tall), for the
efficiency of the HSVD algorithm.  In both cases, the HSVD of a square
matrix will be computed, either of $R$ (in the case of the QR
factorization), or of $R^{\ast}$ (in the case of hyperbolic QR
factorization).  To obtain matrices of singular vectors, after the
completion of the HSVD, either $U$ (in the case of ordinary QR
factorization), or $V$ (in the case of hyperbolic QR factorization)
should be premultiplied by $Q$.

On the other hand, if $M$ is given explicitly, $G$ and $J$ are usually
computed by the Hermitian indefinite factorization with postprocessing
(see \cite{Slapnicar-98}), which always produces $G$ of the full
column rank.  Note that $G$ obtained by the Hermitian indefinite
factorization can be rectangular with more rows than columns, and
should be shortened by the QR factorization.  Therefore, it is
sufficient to efficiently compute HVSD of a square matrix.

The HSVD also serves as the second stage in solving a Hermitian
indefinite eigenproblem.  Our method of choice for such a computation
is the one-sided hyperbolic Jacobi algorithm \cite{Veselic-93}, which
has been much studied in the recent years, due to its high relative
accuracy \cite{Slapnicar-03} and various possibilities of the
efficient blocking and parallelization in the scope of the
conventional CPU computing \cite{Singer-et-al-10,Singer-et-al-10a}.

The efficient GPU solution (especially, one with minimal CPU
intervention) for the first stage (i.e., Hermitian indefinite
factorization) remains an open problem.

We will show that the Jacobi algorithm is elegantly parallelizable on
the modern GPUs, utilizing them as a primary target for the algorithm
execution.  The CPU acts only as the ancillary unit: for driving and
synchronizing the GPU computation, and for data initialization on the
GPU.

The numerical experiments demonstrate the benefits of such an
approach, compared to the fastest existing sequential and
multi-process parallel Jacobi implementations.  Significant speedup
vs.\ the sequential, and moderate speedup vs.\ the CPU-parallel
algorithm with $4$ processes are obtained.

Finally, the true strength of our approach is discussed, through the
possibilities of combining the GPU parallelism with the traditional
one, for the very large problems.

The rest of the paper is organized as follows.  In Section \ref{sec:2}
the essentials of the Jacobi HSVD and basic tools for its
parallelization are presented.  Section \ref{sec:3} deals with the
properties and constraints of a class of GPU computing platforms our
algorithm, named GPUJACH1, is designed for.  The algorithm itself is
detailed in Section \ref{sec:4}, and it is compared, by numerical
testing, with a sequential and a CPU-parallel implementation of the
Jacobi HSVD in Section \ref{sec:5}.  We conclude with a discussion of
the possible applications of GPUJACH1 in the context of hybrid
CPU--GPU parallelism, and with notes on future work, in Section
\ref{sec:6}.  The Appendix contains a proof of convergence for the
modulus pivot strategy.

%
%%%%%%%%%%%%%%%%%%%%%%%%%%%%%%%%%%%%%%%%%%%%%%%%%%%%%%%%%%%%%%%%%%%%%%%%%%%%%%
%
\section{The essentials of the Jacobi HSVD}
\label{sec:2}
%
%%%%%%%%%%%%%%%%%%%%%%%%%%%%%%%%%%%%%%%%%%%%%%%%%%%%%%%%%%%%%%%%%%%%%%%%%%%%%%
%
The hyperbolic one-sided Jacobi algorithm implicitly diagonalizes a
definite pair $(A, J)$, where $A \assgn G^{\ast} G$, by $J$--unitary
congruences \cite{Slapnicar-03}.  ``Implicitly'' in this context means
that the columns of $G$ are orthogonalized.

The idea of the one-sided algorithm is to multiply the columns of
$G^{(0)} := G$ from the right-hand side by $J$-unitary matrices
$V_k^{-\ast}$, $k = 1, 2, \ldots, s$, to obtain the matrix
$G_{}^{(s)}$, with numerically orthogonal (not orthonormal) columns.
If $[V^{-\ast}]^{(0)}$ is defined as $[V^{-\ast}]^{(0)} = I$, then
\begin{equation}
  G^{(k)} = G^{(k - 1)} V_k^{-\ast}, \quad
  [V^{-\ast}]^{(k)} = [V^{-\ast}]^{(k - 1)} V_k^{-\ast}, \quad
  k = 1, 2, \ldots, s.
\label{2.1}
\end{equation}
If $G^{(s)}$ has sufficiently orthogonal columns, $[V^{-\ast}]^{(s)}$
serves as a quite good approximation of $V^{-\ast}$.  Now, the matrix
$V$ of right (hyperbolic) singular vectors is easily obtainable from
$V^{-\ast}$, since $V$ is a $J$-unitary matrix, i.e., $V = J V^{-\ast} J$.

The hyperbolic singular values $\sigma_i$ are norms of the final
$G_{}^{(s)}$, while the approximate matrix $U$ of left singular
vectors is computed by column scaling of $G_{}^{(s)}$ by diagonal
matrix $\Sigma^{-1} = \diag(\sigma_1^{-1}, \ldots, \sigma_n^{-1})$.

If the algorithm is used only for the eigenvalue computation,
$V^{-\ast}$ need not be accumulated, since $U$ keeps the eigenvectors
of $M = G J G^{\ast}$, while the eigenvalues are squares of singular
values multiplied by a correct sign from $J$, i.e.,
$\lambda_i^{} = \sigma_i^2 j_{ii}$.

Just for simplicity, from now on, we restrict ourselves to the real
case, i.e., to symmetric matrices instead of Hermitian.  To emphasize
this, we will use symbol ${}^T$ (instead of ${}^{\ast}$) for
transposition.

The $J$-unitary matrices $V_k^{-T}$ from (\ref{2.1}) are usually
very simple -- they are plane rotations (trigonometric and
hyperbolic).  Such a rotation orthogonalizes a pair of columns of
$G$.  Their non-identity parts can be represented as
\begin{equation}
  \begin{aligned}
    V_T^{-T} & \assgn \begin{bmatrix}
      \hpm \cos \varphi & \sin \varphi \\
          -\sin \varphi & \cos \varphi
    \end{bmatrix}
    = \begin{bmatrix}
      1             & \tan \varphi \\
      -\tan \varphi & 1
    \end{bmatrix}
    \begin{bmatrix}
      \cos \varphi & 0 \\
      0            & \cos \varphi
    \end{bmatrix}, \\[3pt]
    V_H^{-T} & \assgn \begin{bmatrix}
      \cosh \varphi & \sinh \varphi \\
      \sinh \varphi & \cosh \varphi
    \end{bmatrix}
    = \begin{bmatrix}
      1             & \tanh \varphi \\
      \tanh \varphi & 1
    \end{bmatrix}
    \begin{bmatrix}
      \cosh \varphi & 0  \\
      0             & \cosh \varphi
    \end{bmatrix}.
\end{aligned}
\label{2.2}
\end{equation}
Note that $\varphi$ is not needed explicitly, the tangent would
suffice to construct a rotation.  In the Hermitian case, only one
additional angle $\alpha$ is needed, see for example
\cite{SingerSanja-93}.

The column orthogonalization starts by forming a $2 \times 2$ pivot
block $A_p$,
\begin{displaymath}
  A_P = \begin{bmatrix}
    a_{ii} & a_{ij} \\
    a_{ij} & a_{jj}
  \end{bmatrix}
  \assgn \begin{bmatrix}
    g_i & g_j
  \end{bmatrix}^T
  \begin{bmatrix}
    g_i & g_j
  \end{bmatrix},
\end{displaymath}
where $g_i$ and $g_j$ denote the $i$-th and $j$-th column of $G$ (we
may assume $i<j$).

Optionally, a further processing of the column pair may be avoided if
the columns are relatively orthogonal (up to the machine precision
$\varepsilon$), i.e.,
\begin{equation}
  |a_{ij}| < \varepsilon\sqrt{a_{ii} a_{jj}},
\label{2.3}
\end{equation}
where $\varepsilon$ is the smallest floating-point number such that
$\fl(1+\varepsilon) > 1$.

First, the adequate rotation (trigonometric $V_T^{-T}$ if $i$-th and
$j$-th diagonal element of $J$ agree, or hyperbolic $V_H^{-T}$
otherwise) to annihilate the element $a_{ij}$ is computed.  Then, this
transformation (\ref{2.2}) is applied from the right to the columns
$[g_i\ g_j]$, giving
\begin{equation}
    g'_i \assgn (g_i^{} - \tan \varphi g_j^{}) \cos \varphi, \qquad
    g'_j \assgn (g_j^{} + \tan \varphi g_i^{}) \cos \varphi
\label{2.4}
\end{equation}
in the trigonometric case, and
\begin{equation}
    g'_i \assgn (g_i^{} + \tanh \varphi g_j^{}) \cosh \varphi, \qquad
    g'_j \assgn (g_j^{} + \tanh \varphi g_i^{}) \cosh \varphi
\label{2.5}
\end{equation}
in the hyperbolic case.  The formul\ae\ (\ref{2.4})--(\ref{2.5}),
written pointwise, represent a basic tool for updating the columns of
$G$.  The columns of $V^{-T}$ are updated in the same fashion, with an
appropriate change of notation, see (\ref{2.1}).

The parts of the formul\ae\ (\ref{2.4})--(\ref{2.5}) written in parentheses
can be performed by a single fused multiply-add (FMA) operation, where
available, with only one rounding of the result, thus improving the speed,
while exhibiting smaller roundoff errors.

All Jacobi-type algorithms iterate until convergence criteria are
met.  Usually, these iterations are organized in sweeps (sometimes
called cycles).  For a symmetric matrix $A$ of order $n$, a sweep is a
collection of $n (n - 1) / 2$ annihilations of different elements in
the strict upper triangle of $A$.  In other words, in a sweep, pairs
of columns of $G$ (pivot pairs) are orthogonalized as described
previously.  The order of pivot pairs is chosen according to a pivot
strategy.  In sequential use, widespread pivot strategies are row and
column cyclic, since they use cache memory well (depending on layout
of arrays in memory).  For both strategies, convergence
\cite{Veselic-93} and asymptotic quadratic convergence
\cite{Drmac-Hari-93} have been proved.  If the pivot pairs are
orthogonalized in a cyclic manner more than once (but a prescribed
number of times) in a `sweep', this `sweep' is usually called a
quasi-sweep and such a pivot strategy is called quasi-cyclic.

A simple convergence criterion is that no rotations occurred in a
(quasi-)sweep, due to condition (\ref{2.3}), which guarantees relative
accuracy \cite{Slapnicar-03}.  The process could also be stopped if
the computed tangents are all below a predefined threshold
$\Teps := \sqrt{\varepsilon} / 2$, i.e., when the quadratic
convergence is detected.  The second, threshold test is a replacement
for relative orthogonality criterion (\ref{2.3}) from
\cite{Slapnicar-03} since we are trying to avoid an extra
(quasi-)sweep.  If we define
\begin{displaymath}
  \hyp = \begin{cases}
    -1, & \text{in the trigonometric case,} \\
    \hpm 1, & \text{in the hyperbolic case,}
  \end{cases}
\end{displaymath}
then
\begin{displaymath}
  \theta = \begin{cases}
    \tan 2\varphi & \text{if $\hyp = -1$,} \\
    \tanh 2\varphi & \text{if $\hyp = 1$,}
  \end{cases}
  \qquad
  t = \begin{cases}
    \tan \varphi & \text{if $\hyp = -1$,} \\
    \tanh \varphi & \text{if $\hyp = 1$.}
  \end{cases}
\end{displaymath}
If the tangent is very small, i.e., $|t| \ll 1$, then
$\theta \approx 2t$.  The angle $\theta$ is computed from the
requirement that $a'_{ij} = 0$,
\begin{displaymath}
  \theta = \frac{2 a_{ij}}{-a_{ii} + \hyp \cdot a_{jj}}.
\end{displaymath}
Therefore, for $|t| \leq \Teps$, we have
$a_{ij} \approx t (-a_{ii} + \hyp \cdot a_{jj})$, and since
\begin{equation}
    a'_{ii} = a_{ii}^{} + \hyp \cdot t a_{ij}^{}, \qquad
    a'_{jj} = a_{jj}^{} + t a_{ij}^{},
\label{2.6}
\end{equation}
when $|t| \leq \Teps$, the second term from the right-hand side in
(\ref{2.6}) satisfies
\begin{displaymath}
  |\hyp \cdot t a_{ij}| =  |t a_{ij}| \lesssim \frac{\varepsilon}{4} \cdot
  (|-a_{ii} + \hyp \cdot a_{jj}|)
  \leq \frac{\varepsilon}{2} \cdot \max \{a_{ii}, a_{jj} \}.
\end{displaymath}
Therefore $a'_{ii} \approx a_{ii}^{}$ and $a'_{jj} \approx a_{jj}^{}$
in (\ref{2.6}).  When $\diag(JA)$ is sorted and the process is
near the end, the off-diagonal norm of $A$, i.e., $\|A-\diag(A)\|_F$,
quadratically diminishes \cite{Drmac-Hari-93}, and so the off-diagonal
elements.  The tangents in the following \hbox{(quasi-)}sweep will also
be quadratically smaller, and the diagonal updates in (\ref{2.6})
will be negligible.

To summarize, the orthogonality convergence check is a safe fallback,
with guaranteed numerical orthogonality of the eigenvectors, but in our
testing experience, the threshold criterion has always happened before
and terminated the process.

The Jacobi HSVD algorithm, when computing the eigendecomposition, can
be implemented using just one 2-dimensional array \texttt{G},
initialized to the factor $G$.  At the end of the process, \texttt{G}
holds the computed $U\Sigma$.  The column norms of $U\Sigma$ are the
hyperbolic singular values, and normalized columns are the left
singular vectors of $G$, i.e., the eigenvectors of $M$.  If the full
HSVD is desired, the sequence of applied $V_T^{-T}$ and $V_H^{-T}$
transformations, i.e., the matrix $V^{-T}$, needs to be accumulated in
another 2-dimensional array \texttt{V}, starting with the identity
matrix $I$.

On a theoretical level, the main difference between a sequential and
a parallel one-sided Jacobi algorithm is the choice of a pivot
strategy.  In any parallel algorithm, CPU or GPU based, a pivot
strategy is chosen to simultaneously orthogonalize as many independent
(non-overlapping) pivot pairs as possible.  These pairs can be either
two single columns, or two block-columns.  Either way, this collection
of independent pairs will be called a step.

Our choice of parallel pivot strategy, for CPU and GPU algorithms
alike, is a slightly modified modulus strategy \cite{LukF-Park-89}.
If the order $n$ of a matrix is even, $n=2q$, then a quasi-sweep has
exactly $n$ steps, and all the elements in the upper triangle of $A$
are annihilated (some of them twice).  For $q$ independent tasks
(e.g., CPU or GPU threads), exactly $q$ pivot pairs are simultaneously
orthogonalized in each step.

\begin{figure}[hbt]
\begin{center}
  \includegraphics{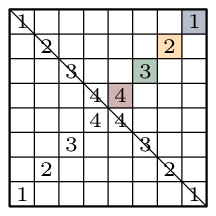}\hfil
  \includegraphics{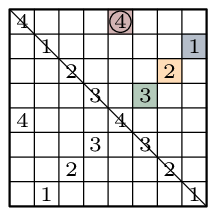}\hfil
  \includegraphics{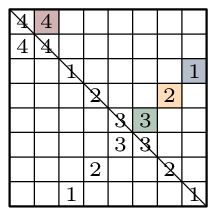}\hfil
  \includegraphics{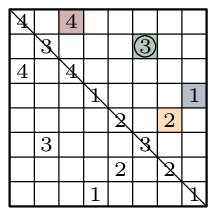}\hfil
  \phantom{\includegraphics{fig4.eps}}\\[6pt]
  \includegraphics{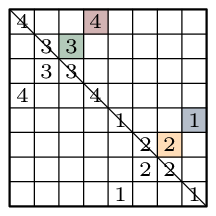}\hfil
  \includegraphics{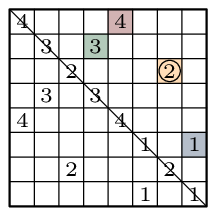}\hfil
  \includegraphics{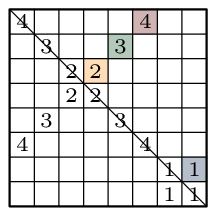}\hfil
  \includegraphics{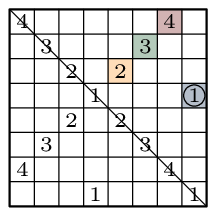}\hfil
  \includegraphics{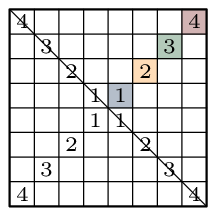}\hfil
\end{center}
\caption{Modified modulus strategy for $n=8$ and $q=4$.
  Background colors denote elements annihilated in one step,
  and circled elements are annihilated twice in a quasi-sweep.}
\label{fig:2.1}
\end{figure}

In Fig.~\ref{fig:2.1} the modified modulus pivot strategy, for $n=8$
($q=4$), is illustrated from a point of view of the (implicit) matrix
$A$.  Each subfigure shows a single step, with $4$ pivot submatrices
that can be independently processed by $4$ tasks, either sequentially
or simultaneously.  After $8$ steps a quasi-sweep is completed, and a
new one will begin with the initial assignment of pivot columns to
tasks reversed, as shown in the ninth subfigure.

Our modified modulus pivot strategy differs from \cite{LukF-Park-89}
in the choice of a start step (ours is the antidiagonal), and in
annihilating exactly twice the elements of the diagonal $(i,q+i)$,
$i=1,\ldots,q$, for even $n$.  This is a quasi-cyclic strategy,
designed to speedup the convergence.  The proof of convergence for
this strategy is a work in progress.  If we avoid double annihilation,
our strategy is shift equivalent to the modulus strategy, for which
the proof of convergence is detailed in Appendix.

For more details of the modified modulus pivot strategy, including the
efficient algorithm for step transitions, please see Section
\ref{sec:4} and \cite{Singer-et-al-10a}.

%
%%%%%%%%%%%%%%%%%%%%%%%%%%%%%%%%%%%%%%%%%%%%%%%%%%%%%%%%%%%%%%%%%%%%%%%%%%%%%%
%
\section{The GPU computing platform}
\label{sec:3}
%
%%%%%%%%%%%%%%%%%%%%%%%%%%%%%%%%%%%%%%%%%%%%%%%%%%%%%%%%%%%%%%%%%%%%%%%%%%%%%%
%
GPU computing has already evolved to a mainstream technology.
Although vendor-neutral standards (like OpenCL \cite{Khronos-10}) have
emerged, we implemented GPUJACH1 for the NVIDIA CUDA architecture
\cite{NVIDIA-10}, due to its maturity and close ties to the underlying
hardware.  The algorithm is portable, more or less efficiently, to any
similar single-instruction multiple-threads platform (SIMT), provided
IEEE 754-2008 \cite{IEEE-08} double-precision floating-point
arithmetic is available and the basic computational concepts are found
in (or could be mapped to) the other architecture.

In short, the CUDA architecture provides to a programmer a thin layer
of abstractions, programming tools and interfaces on top of the recent
NVIDIA graphics processing hardware.  The hardware itself is seen as a
massively parallel device whose threads execute the same instruction
sequence over (possibly) different data, stored on the device.  The
execution is initiated by the CPU (host), and a fast host-to-device
and device-to-host data transfer is available (but not as fast as
device memory access).

This SIMT paradigm entails a careful rethinking of even the basic
algorithms.  Any branching causes a significant slowdown, since
threads on divergent branches proceed sequentially.  The Jacobi
algorithm, however, fits this paradigm perfectly, because the same
orthogonalization primitives are applied to different pairs of matrix
columns concurrently.

A subroutine the device threads are executing is called a kernel.  A
kernel can be written in a high-level programming language (e.g., C),
or using the assembly-like, low-level, but general-purpose instruction
set of the PTX \cite{NVIDIA-10a} virtual machine.

The device threads are grouped into so-called warps, for memory access
optimization, execution scheduling and finest-grain synchronization.
A warp consists of $32$ threads in the current CUDA hardware, with two
half-warps of $16$ threads each.

From a high-level perspective, blocks of threads are established.  All
blocks are of the same size, and may be seen as one, two or three
dimensional arrays of threads.  Threads are indexed and may be
synchronized only within their blocks, i.e., different blocks are
concurrent and mutually independent.  For each launch of a kernel, the
size of a block and the number of blocks are set, tailored to the
nature of the algorithm and the data.  Such a one or two dimensional
array of blocks is called a grid.

\begin{figure}[hbt]
\begin{center}
  \includegraphics{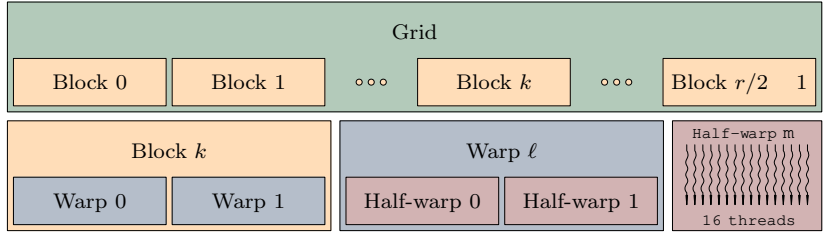}
\end{center}
\caption{The CUDA grid for GPUJACH1.}
\label{fig:3.1}
\end{figure}

The grid for GPUJACH1 kernels, working on a factor $G$ with even
number of columns $r$, is shown in Fig.~\ref{fig:3.1}.

Memory management is the major obstacle in GPU programming, as many
distinct memory spaces of varying speed, size and accessibility are
involved.

Input and output data is stored in the global memory, a large but slow
portion of the GPU memory.  It is accessible to the CPU and all GPU
threads.  Improper (uncoalesced) access tremendously degrades
performance \cite{NVIDIA-10}, thus the carefully aligned addressing of
successive locations by threads with successive indices is required.

The shared memory, on the other hand, is a small but fast device
memory, allocated per block, and ideal for data exchange between
threads of the same block.

The constant memory, small and cached, holds kernel parameters not
changing between launches (e.g., device pointers to global memory
arrays in GPUJACH1).

Arithmetic is done in per-thread registers.  GPUJACH1 uses 64-bit
floating-point and 32-bit integer arithmetic.  FMA instruction (e.g.,
dot-products, column updates) and 24-bit integer multiplication (e.g.,
address calculations) are chosen, if possible.

GPUJACH1 targets NVIDIA GT200 graphics processor series, common in the
time of writing.  The threads are executed by $30$ 8-core
multiprocessors (SMs).  Each SM has a $64\rm\ kB$ register file, one
double-precision arithmetic unit and $16\rm\ kB$ of shared memory.
More than one block can reside on an SM, as the resources (register
usage per thread and shared memory allocation per block) and other
constraints \cite{NVIDIA-10} allow.  GPUJACH1 needs only $512\rm\ B$
of shared memory per block, but $48$ 32-bit registers per thread, thus
at most $5$ blocks can reside on an SM at a given time.

Two distinctive features of the GT200 series are that no cache is
available for the global memory, and that single-precision arithmetic
deviates in the guaranteed accuracy too much from the standard
\cite{IEEE-08} to be useful for GPUJACH1.  Both issues are addressed
by the more advanced architectures (e.g., NVIDIA Fermi
\cite{NVIDIA-10}).

%
%%%%%%%%%%%%%%%%%%%%%%%%%%%%%%%%%%%%%%%%%%%%%%%%%%%%%%%%%%%%%%%%%%%%%%%%%%%%%%
%
\section{The GPUJACH1 algorithm}
\label{sec:4}
%
%%%%%%%%%%%%%%%%%%%%%%%%%%%%%%%%%%%%%%%%%%%%%%%%%%%%%%%%%%%%%%%%%%%%%%%%%%%%%%
%
The GPUJACH1 algorithm is an efficient parallel realization of the
one-sided pointwise hyperbolic Jacobi SVD.  GPUJACH1 provides all
parts of the Jacobi process executed on the GPU as much as possible.

GPUJACH1 consists of, essentially, two parts:
\begin{compactitem}
\item\emph{host\/} routines (auxiliary), which are executed on the CPU, and
\item\emph{device\/} routines (computational), which are executed on the GPU.
\end{compactitem}
Host routines serve mainly to call device routines in a synchronized
manner and collect the resulting information, where needed.

A sketch of the main driver routine is given in Alg.~\ref{alg:4.1},
and then its essential parts are described in details.  In the
following, by the subscript $H$ host variables, and by $D$ device
variables are denoted.  Arrays are assumed to be in Fortran
(column-major) order, but indices are zero-based (as in C).  Array
slices are written in Matlab fashion.

\begin{algorithm}
\SetKwInput{KwInOut}{Input/Output}
\SetKwInput{KwDesc}{Description}
\KwDesc{This is the main program, executed on the CPU.}
\KwIn{$n\times r$ factor $G$, and the number $p$ of positive signs in $J$.}
\KwOut{hyperbolic singular values $\Sigma$ and the matrices of left ($U$) and
  hyperbolic right ($V^T$) singular vectors.}
\SetKwFunction{HostDriverRoutine}{Host\_Driver\_Routine}
\SetKwFunction{PrecomputeData}{Precompute\_Data}
\SetKwFunction{JacobiStep}{Jacobi\_Step}
\SetKwFunction{SortDiagonal}{Sort\_Diagonal}
\SetKwFunction{ResetConvergence}{Reset\_Convergence}
\SetKwFunction{CheckConvergence}{Check\_Convergence}
\SetKwFunction{CollectDataToHost}{Collect\_Data\_To\_Host}
\SetKwFunction{LoadDataToDevice}{Load\_Data\_To\_Device}
\SetKwComment{Comment}{\quad /\!/\,}{}
\BlankLine
\HostDriverRoutine\;
\Begin{
    $G_D \leftarrow G_H$; \quad $V_D \leftarrow I_r$\Comment*[l]{optional}
    \PrecomputeData\;
    \SortDiagonal\;
    \For{$\textit{sweep\/} = 0 \ \KwTo\  \textit{MaxSweep\/} - 1$}
        {
          \ResetConvergence\;
          \For{$\textit{step\/} = 0 \ \KwTo \ \textit{r\/} - 1$}
              {
                \JacobiStep\;
              }
          \CheckConvergence\;
          \SortDiagonal\;
        }
    $G_H \leftarrow G_D$; \quad $V_H \leftarrow V_D$; \quad $D_H \leftarrow D_D$\Comment*[l]{optional}
  }
\caption{The host driver routine}
\label{alg:4.1}
\end{algorithm}

GPUJACH1 holds the $n\times r$ factor $G$ in arrays $G_{H,D}$ (we
write $G_{H,D}$ as a shorthand for $G_H$ and $G_D$, similarly for
$V$).  We may assume the factor to be square ($r\times r$), as left by
preprocessing (see Section \ref{sec:1}).  The factor could also be
preloaded on the GPU by a previous processing, and not needed on the
CPU afterward, so $G_H$ is optional.  The factor must be of the full
column rank, $r$ needs to be even and not greater than $n$, and $n$,
for performance reasons, should be a multiple of the warp size.  This
is not a loss of generality, since the initial $G$ could easily be
bordered, as in (\ref{4.1}), to satisfy these constraints
\begin{equation}
  G'_{H, D} =
  \begin{array}{c}
    n \, \left\{
    \begin{array}{r}
      \hskip0ex \vrule width 0cm height 0.59cm depth 0.59cm
    \end{array}
    \right. \\
    \hskip-1ex \vrule width 0cm height 0.64cm depth 0.64cm
  \end{array}
  \hskip-1ex
  \left[ \begin{array}{c:c}
      \begin{array}{ccc}
        \phantom{G} &   & \phantom{G} \\
	    \phantom{G} & G & \phantom{G} \\
        \phantom{G} &   & \phantom{G}
      \end{array}
      & \hphantom{\,\,} 0 \\
      \hdashline
      \begin{array}{ccc}
        \phantom{G}\vphantom{\big|} & \strut 0 & \phantom{G}
      \end{array}
      & \hphantom{\,\,}\vphantom{\big|}\strut 1 \\
      \hdashline
      \begin{array}{ccc}
        \phantom{G} & \vphantom{\bigg|} 0 & \phantom{G}
      \end{array}
      &
      \begin{array}{c}
        \hphantom{\,\,} \vphantom{\bigg|} 0
      \end{array}
    \end{array} \right]
  \left.
  \begin{array}{l}
    \hskip-1ex \vrule width 0cm height 1.3cm depth 1.3cm
  \end{array}
  \right\} \ \ell \cdot \textit{WarpSize\/}.
  \label{4.1}
\end{equation}

If the matrix $V^{-T}$ is needed, it is accumulated in the array
$V_D$, starting from the identity $I_r$, and is optionally transferred
to $V_H$ at the end.  Arrays $V_{H,D}$, if used, must have the same
number of columns as $G_{H,D}$, and the same restrictions (and
appropriate bordering) apply.

In principle, the execution begins with all the data residing in the
main (CPU) memory, and should be copied to the global memory of the
GPU.

\begin{figure}[hbt]
\begin{center}
  \includegraphics{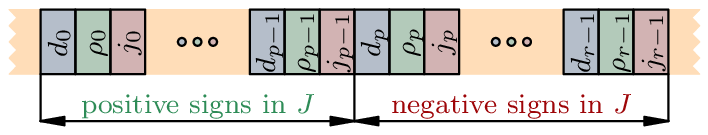}
\end{center}
\caption{$D_{H,D}$ -- packing of diagonal, permutation, and sign
  $r$-vectors in one.}
\label{fig:4.1}
\end{figure}

For efficiency, to access global memory data as few times as possible,
and to reuse results of the computation while still in registers, the
diagonal of the implicit pivot matrix $A$ is kept and updated in each
Jacobi step.  While performance gains are here obvious, yet additional
speedup can be obtained by a special diagonal sorting, described in
\cite{Singer-et-al-10} for block versions of the Jacobi algorithm.
Keeping the diagonal, i.e., the eigenvalues, sorted ensures the
quadratic convergence of the Jacobi algorithm
\cite{Drmac-Hari-93,Singer-et-al-07}.  To facilitate the sorting, the
diagonal entries $d$, the composition of the sorting permutations
$\rho$, and the signature $J$ is packed into vectors $D_{H,D}$, as
shown in Fig.~\ref{fig:4.1}.

The separation of $D_{H,D}$ to $p$ entries with positive and $r-p$
entries with negative sign elements of $J$ is always maintained.
$D_D$ is precomputed in Alg.~\ref{alg:4.2}, and sorted in
\texttt{Sort\_Diagonal} routine.  The sorting key of a diagonal
package is $d$.  The first part of the diagonal (the one with positive
signs in $J$) is sorted decreasingly, and the last part (with negative
signs in $J$) is sorted increasingly, with respect to that key.  We
chose a GPU implementation of the merge-sort, \texttt{thrust::sort},
from the Thrust library \cite{Thrust}, to serve as
\texttt{Sort\_Diagonal}.  Note that, at the end of the process, $D_D$
holds the sorted and squared hyperbolic singular values $d$ and the
final permutation $\rho$.

\begin{algorithm}
\SetKwInput{KwInOut}{Input/Output}
\SetKwInput{KwDesc}{Description}
\KwDesc{Routine computes
  $D_D^{} = (d=\diag(G_D^TG_D^{}),\ \rho=\id_r,\ j=\diag(I_p,-I_{r-p}))$,
  and initializes the stepper vector $S_D$.}
\SetKwFunction{PrecomputeData}{Precompute\_Data}
\SetKwFunction{Device}{Device\_Code}
\SetKwFunction{InitStepper}{Initialize\_Stepper}
\SetKwComment{Comment}{\quad /\!/\,}{}
\SetKwComment{xComment}{\ /\!/\,}{}
\BlankLine
\PrecomputeData(\Device(block $k$))\;
\Begin{
    $(i, j) = (2k, 2k + 1)$\;
    $d_i^{} = g_i^T g_i^{}$; \quad $d_j^{} = g_j^T g_j^{}$\;
    $\rho_i = i$; \quad $\rho_j = j$\;
    \lIf{$i < p$}{$j_i = 1$}
    \lElse{$j_i = -1$}; \quad
    \lIf{$j < p$}{$j_j = 1$}
    \lElse{$j_j = -1$}\;
    \InitStepper\;
  }
\caption{Precompute\_Data routine}
\label{alg:4.2}
\end{algorithm}

Diagonal packing from Fig.~\ref{fig:4.1} might seem unnecessary, when
only diagonal entries and the permutation are strictly needed, and
could be arranged into two simple, separate vectors of doubles and
integers, respectively.  Diagonal entries could also carry the
respective signs from $J$, to simplify sorting.  This approach is
certainly possible, but it shows in the following section to be
slower, due to slightly more complex \texttt{Jacobi\_Step} routine,
than the presented solution.  Moreover, packed $j$ incurs no extra
overhead over packing solely $d$ and $\rho$ together -- the place $j$
occupies in the package must exist because of data alignment
constraints of the GPU \cite{NVIDIA-10}.

The CUDA grid for GPUJACH1 kernels (Fig.~\ref{fig:3.1}) has exactly
$b = r / 2$ blocks.  Each block is assigned its own pair of columns of
$G_D$.  In each Jacobi step these pivot pairs are mapped to blocks
according to the modified modulus pivot strategy, implemented as an
integer vector
\begin{displaymath}
  S_D = (\ip, \jp, \iblk, \jblk)
\end{displaymath}
%$S_D=(\ip,\jp,\iblk,\jblk)$
of length $4 \times b$, and
a device routine, called \texttt{Update\_Stepper}, which updates $S_D$
at the end of a step, independently by each block (see
Algs.~\ref{alg:4.3} and \ref{alg:4.4}).

\begin{algorithm}
\SetKwInput{KwDesc}{Description}
\SetKwFunction{Device}{Device\_Code}
\SetKwFunction{InitStepper}{Initialize\_Stepper}
\InitStepper(\Device(block $k$))($r$, $S_D=(\ip,\jp,\iblk,\jblk)$)\;
\KwDesc{Initialization of the modified modulus strategy stepper vector
  $S_D$ to antidiagonal -- $k \leftrightarrow (k, r - k - 1)$.
  Variables $\ip$ and $\jp$ are auxiliary, while $\iblk_k$ and $\jblk_k$
  contain the first and the second column index (non-permuted) of a
  pivot pair for block $k$.}
\BlankLine
\Begin{
    $\ip_k = k$;\quad
    $\iblk_k = \ip_k$;\quad
    $\jp_k = r - k - 1$;\quad
    $\jblk_k = \jp_k$\;
  }
\caption{Initialize\_Stepper routine}
\label{alg:4.3}
\end{algorithm}

\begin{algorithm}
\SetKwComment{Comment}{\quad /\!/\,}{}
\SetKwFunction{Device}{Device\_Code}
\SetKwFunction{UpdateStepper}{Update\_Stepper}
\UpdateStepper(\Device(block $k$))($r$, $S_D=(\ip,\jp,\iblk,\jblk)$)\;
\SetKwInput{KwDesc}{Description}
\KwDesc{Block $k$ determines the indices of the next pivot pair
  $(\iblk_k, \jblk_k)$.}
\BlankLine
\Begin{
    \uIf
        {$(\ip_k + \jp_k) > r - 1$}
        {$\ip_k = \ip_k + 1$\;
          \If
             {$\ip_k = \jp_k$}
             {$\ip_k = \ip_k - r/2$;\quad $\jp_k = \ip_k$\;}
         $\iblk_k = \ip_k$\;
        }
    \Else
        {$\jp_k = \jp_k + 1$;\quad $\jblk_k = \jp_k$\;}
  }
\caption{Update\_Stepper routine}
\label{alg:4.4}
\end{algorithm}

In other words, each block orthogonalizes its pivot pair
$(\rho(\iblk_k),\rho(\jblk_k))$.  In a block, each warp is dedicated
to access a single column only.  This design makes the device code
almost completely uniform for all threads, thus avoiding branching as
much as possible.

In the \texttt{Jacobi\_Step} (Alg.~\ref{alg:4.5}) columns of $G_D$ and
$V_D$ are indexed by the current permutation $\rho$, and no physical
swapping of columns ever takes place.  Therefore, at the end of the
process, $G_D$ holds $U\Sigma$, and $V_D$ holds $V^{-T}$, in the
original column order.  To match the computed (and permuted) singular
values to the singular vectors, $D_{H,D}$ need to be permuted by the
inverse of the final permutation $\rho$.

\begin{algorithm}
\SetKwInput{KwInOut}{Input/Output}
\SetKwInput{KwDesc}{Description}
\KwDesc{This is the main GPU computational routine.}
\SetKw{KwGoto}{goto}
\SetKwFunction{JacobiStep}{Jacobi\_Step}
\SetKwFunction{UpdateCvg}{Update\_Convergence}
\SetKwFunction{UpdateStepper}{Update\_Stepper}
\SetKwFunction{End}{End}
\SetKwFunction{Device}{Device\_Code}
\SetKwComment{Comment}{\quad /\!/\,}{}
\BlankLine
\JacobiStep(\Device(block $k$))\;
\Begin{
    $(i,j) = (S_D(k).\iblk, S_D(k).\jblk)$\Comment*[l]{$i < j$}
    $a_{ii} = d_i$; \quad $a_{jj} = d_j$\Comment*[l]{already computed}
    $a_{ij}^{} = g_{\rho(i)}^T g_{\rho(j)}^{}$\Comment*[l]{dot-product}
    \Comment*[h]{Indexing by permutation; no physical swapping of $G_D$ columns}\;
    \lIf{
      $g_{\rho(i)}$ and $g_{\rho(j)}$ relatively orthogonal (up to
      given $\varepsilon$)
      }
    {\KwGoto \End}\;
    \lIf{
       $j_i = j_j$}
    {set $\hyp = -1$}
    \lElseIf{
      $j_i = -j_j$}
    {set $\hyp = 1$}\;
    \lIf{
      $\hyp = -1$}
    {compute $t$ from $a_{ii}$, $a_{jj}$, $a_{ij}$ as $\tan \varphi$}
    \lElseIf{
      $\hyp = 1$}
    {as $\tanh \varphi$}\;
    \lIf{
      $\hyp = -1$}
    {compute $c$ from $t$ as $\cos \varphi$}
    \lElseIf{
      $\hyp = 1$}
    {as $\cosh \varphi$}\;
    $g'_{\rho(i)} = (g_{\rho(i)}^{} + \hyp\cdot t g_{\rho(j)}^{}) c$; \quad
    $g'_{\rho(j)} = (t g_{\rho(i)}^{} + g_{\rho(j)}^{}) c$\Comment*[l]{\texttt{FMA}, \texttt{scal}}
    \Comment*[h]{While updating $G_D^{}$ columns, compute new $d'_i$ and  $d'_j$}\;
      $d'_i=(g'_{\rho(i)})^Tg'_{\rho(i)}$; \quad
      $d'_j=(g'_{\rho(j)})^Tg'_{\rho(j)}$\;
    $v'_{\rho(i)} = (v_{\rho(i)}^{} + \hyp\cdot t v_{\rho(j)}^{}) c$; \quad
    $v'_{\rho(j)} = (t v_{\rho(i)}^{} + v_{\rho(j)}^{}) c$\Comment*[l]{\texttt{FMA}, \texttt{scal}}
    \UpdateCvg\Comment*[l]{as in (\ref{4.2})}
    \End: \UpdateStepper\;
    }
\caption{Jacobi\_Step -- the main computational routine}
\label{alg:4.5}
\end{algorithm}

The rest of \texttt{Jacobi\_Step} (Alg.~\ref{alg:4.5}) is more or less
standard.  The final question and the crucial efficiency issue is how
dot products and column updates (\texttt{daxpy}-like operations) are
performed.  We shall describe only the dot product computation, since
updating the columns is done in a similar fashion.

\begin{figure}[hbt]
\begin{center}
  \includegraphics{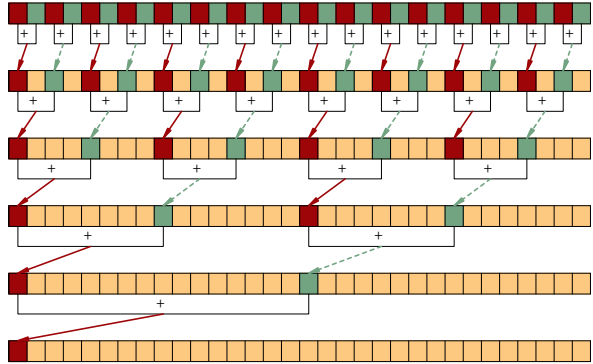}
\end{center}
\caption{Per-warp reduction (via shared memory) of the partial sums to the
  final dot-product.  After each reduction stage ``dashed'' threads
  wait for the whole reduction to complete (but do not terminate).}
\label{fig:4.2}
\end{figure}

For a dot product, each warp ``grabs'' ${\it WarpSize} = 32$
successive elements of its column vector, one element per each thread.
The threads then exchange these values via shared memory and update
(FMA) their local partial sums.  Finally, the partial sums are reduced
in the shared memory, as shown in Fig.~\ref{fig:4.2}.  This reduction
at warp level is free of synchronization \cite{NVIDIA-10}, and needs
to keep all the threads in a warp alive for further processing.

The convergence statistics is held in an integer vector $C_D$ of
length $b$, which is updated independently for each block in a step.
If no rotation was applied in block $k$, $C_D(k)$ remains unchanged,
otherwise it is updated according to the following formula
\begin{equation}
  C_D'(k) = \begin{cases}
    (11)_2            & \text{if $|t| > \Teps$,} \\
    (01)_2 \bor C_D(k)& \text{if $|t| \leq \Teps$,}
  \end{cases}
\label{4.2}
\end{equation}
where $t$ is computed tangent in that block.
Routine \texttt{Reset\_Convergence} zeroes out $C_D$, and
\texttt{Check\_Convergence} returns bitwise-or reduction of $C_D$.
For that reduction on the GPU we used \texttt{thrust::reduce} method
of the Thrust library \cite{Thrust}.

The reduction result tells whether there were no rotations in a
quasi-sweep because the columns are orthogonal up to the machine
precision (i.e., $(00)_2$), or the quadratic convergence was detected
(i.e., $(01)_2$), or we must proceed further since no convergence
criteria were fulfilled (i.e., $(11)_2$).  Otherwise, we stop the
process.

It is worth noting that the convergence status could be accumulated by
atomic bitwise-or global memory instruction \cite{NVIDIA-10a}, but
this approach was significantly slower.

%
%%%%%%%%%%%%%%%%%%%%%%%%%%%%%%%%%%%%%%%%%%%%%%%%%%%%%%%%%%%%%%%%%%%%%%%%%%%%%%
%
\section{Numerical testing}
\label{sec:5}
%
%%%%%%%%%%%%%%%%%%%%%%%%%%%%%%%%%%%%%%%%%%%%%%%%%%%%%%%%%%%%%%%%%%%%%%%%%%%%%%
%
We found it intriguing to hand-code the main GPUJACH1 kernels
(\texttt{Jacobi\_Step} and \texttt{Precompute\_Data}) in the PTX
instruction set and wrap it up with the CUDA Driver API.  The host
part is written in C99 and C++.

GPUJACH1 was compared to the sequential (with row cyclic pivot
strategy) \cite{Slapnicar-03} and MPI-parallel block-oriented
\cite{Singer-et-al-10} one-sided hyperbolic Jacobi algorithms.

The parallel block-oriented algorithm is a blocked CPU version of the
pointwise algorithm described here.  In a step, instead of elements,
(modified) modulus strategy reaches blocks.  In the first step, inside
each block $A_P$ (formed of block columns $G_P \assgn [G_i\ G_j]$ as
$A_P = G_P^{T} G_P^{}$), each element in the strict upper triangle of
$A_P$ is annihilated only once.  In all other steps, only the elements
of the off-diagonal block $G_i^T G_j^{}$ of $A_P$ are annihilated.  In
both cases the inner pivoting strategy is row-cyclic.
To achieve both, matrix $A_p$ is factored by the specially pivoted
Cholesky factorization $P^T A_P^{} P = R^{T} R$, and then a HSVD of
$R$ is computed and the transformation matrix is applied from the
right on $G_P^{}$.  This blocking serves twofold purpose.  The first
is speeding up the process by keeping data (blocks or their parts)
in the cache and utilizing BLAS 3 instead of BLAS 1 operations.
The second is efficiency of the parallelization, since the number of
parallel tasks is usually much smaller than the matrix order.

This approach is the fastest known Jacobi-type parallel algorithm for
computing the HSVD of the moderate-size matrices.  Note that the
pointwise algorithms are less amenable for application on a cluster,
since they require heavy communication between parallel tasks, when
they reside on different machines.

Testing is performed on GTX280, GTX275, and TeslaC1060 graphics cards,
and on two quad-core machines: Intel Xeon X5470 (sequential and 4-task
parallel testing), and Intel Core i7 950 (8-task parallel testing,
with hyperthreading).  For reference purposes, the speed benchmarks
(measured in GFlops) are given in Table \ref{tbl:5.1}.

\begin{table}[hbt]
  \begin{center}
    \begin{tabular}{@{}cccccccccc@{}}
      \toprule
                  &         & Linpack   & \multicolumn{3}{c}{matrices of order $1000 \times 1000$} \\
                  & Machine & Benchmark & \multicolumn{2}{c}{Intel MKL \texttt{dgemm}} & Parallel block-oriented\\
                  & Rpeak   & Rmax      & sequential & parallel                        & hyperbolic Jacobi \\
      \midrule
      Xeon X5470  & 53.328  & 44.581    & 12.267     & 42.535   & 29.951 \\
      Core i7-950 & 48.960  & 45.972    & 12.315     & 39.939   & 35.343 \\
      \bottomrule
    \end{tabular}
  \end{center}
  \caption{Various measures of speed for the testing machines.}
  \label{tbl:5.1}
\end{table}

For testing purposes, symmetric indefinite matrices with random,
uniformly distributed spectra in
$[-a, -a \cdot 10^{-5} \rangle \cup \langle a \cdot 10^{-5}, a]$,
were generated by a modified LAPACK \texttt{dlagsy} routine,
rewritten, with all its support LAPACK and BLAS routines, in 80-bit
\texttt{extended} precision (vs.\ 64-bit \texttt{double}).  Then,
matrices were factorized by the symmetric indefinite factorization
with complete pivoting \cite{Bunch-Parlett-71} (\texttt{xsybpc}), also
in \texttt{extended} precision.  The factor $G$ is then downcasted to
\texttt{double}.  Graphically:
\begin{displaymath}
\stackrel{\texttt{dlarnd}}{\llrightarrow} \; \Lambda_{64} \;
\stackrel{\textrm{upcast}}{\llrightarrow} \; \Lambda_{80} \;
\stackrel{\texttt{xlagsy}}{\llrightarrow} \; M_{80} \;
\stackrel{\texttt{xsybpc}}{\llrightarrow} \; G_{80} \;
\stackrel{\textrm{downcast}}{\llrightarrow} \; G_{64}.
\end{displaymath}
The \texttt{extended} support came from GNU Fortran with
\texttt{real(kind=10}) datatype.  The parameter $a$ varied
depending on the matrix order $n$ as follows in Table \ref{tbl:5.2}.

\begin{table}[hbt]
  \begin{center}
    \begin{tabular}{@{}ccccc@{}}
      \toprule
      $n$ &  $\le 3168$ &  $\le 6368$ &  $\le 9568$ & $\le 10144$\\
      \midrule
      $a$ &        $20$ &        $30$ &        $40$ &        $50$\\
      \bottomrule
    \end{tabular}
  \end{center}
  \caption{Dimension of matrices $n$, and parameter $a$ used in generation.}
  \label{tbl:5.2}
\end{table}

Besides these cases, with approximately the same number of positive
and negative signs in $J$, we also have experimented on matrices with
different number of positive (negative) signs.  The experiments
confirm that the former cases are the slowest ones, while the definite
cases (with zero, or all, signs positive) are the fastest.  An example
of this behavior is given in Table~\ref{tbl:5.3}.

\begin{table}[hbt]
  \begin{center}
    \begin{tabular}{@{}cccccccccc@{}}
      \toprule
      number of positive signs &   0 &   1 &   5 &  10 &  50 & 100 & 500 & 1000 & 2048 \\
      \midrule
      time         (with sorting)    & 273 & 274 & 276 & 285 & 312 & 321 & 354 &  356 &  357 \\
      quasi-sweeps (with sorting)    &  13 &  13 &  13 &  14 &  15 &  15 &  16 &   16 &   16 \\
      \midrule
      time         (without sorting) & 306 & 307 & 308 & 321 & 333 & 349 & 360 &  372 &  368 \\
      quasi-sweeps (without sorting) &  13 &  13 &  13 &  14 &  14 &  15 &  15 &   16 &   16 \\
      \bottomrule
    \end{tabular}
  \end{center}
  \caption{Speed of HSVD on GTX275 for matrices of order $4096$
    (with accumulation of $V^{-T}$) and different number of
    positive signs in $J$.}
  \label{tbl:5.3}
\end{table}
\vskip 0pt plus0pt minus3pt

This behavior is well-known \cite{Hari-SingerSanja-SingerSasa-10}.
The reason lies in fact that trigonometric transformations keep the
trace of the matrix constant, while the hyperbolic transformations
lower the trace.  In exchange, the hyperbolic transformations enables
the high relative accuracy of computed eigenvalues, and makes the
two-norm of the spectral projector small when the eigenvalues have big
relative gaps (see \cite{Slapnicar-Veselic-95}).

The spectrum $\Lambda$ was saved in all cases, to be compared with
the computed eigenvalues $\Lambda'$.  The orders of matrices were
$160 + 128k$, with $0\le k \le 78$.  Large values of $k$ were tested
only on TeslaC1060 (with 4 GB of memory), and skipped on GTX275 and
GTX280, due to the cards' insufficient memory sizes of 896 MB and
1 GB, respectively.

\begin{figure}[hbt]
\begin{center}
  \includegraphics{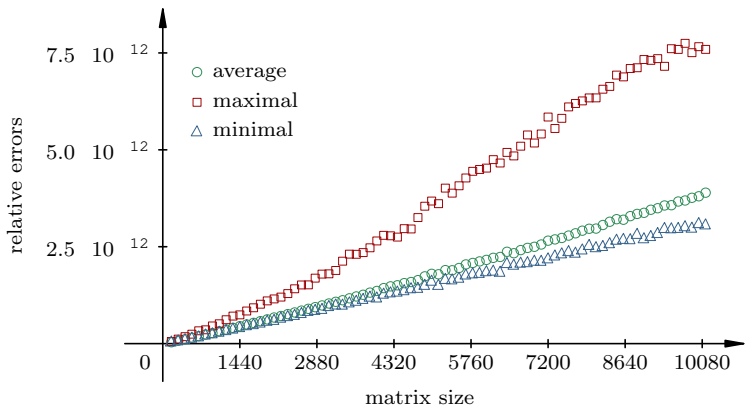}
\end{center}
\caption{Relative errors in the computed eigenvalues,
  $\displaystyle\frac{|\lambda_i^{}-\lambda'_i|}{|\lambda_i^{}|}$.}
\label{fig:5.1}
\end{figure}

Relative errors in the computed eigenvalues (Fig.~\ref{fig:5.1}) are
satisfactory in the context of high relative accuracy.  Note that the
average error is close to minimal.

Distance from orthonormality for the computed eigenvectors $U$ is
calculated as $d(U)=\|I - U^T U\|_F$ (see \cite{Mathias-96} for
details).  It grows linearly form $1.11 \cdot 10^{-14}$ to
\hbox{$7.55 \cdot 10^{-13}$} on full test spectrum.
Row-distance $d(U^T)$ differs from $d(U)$, consistently with
the bound $d(U)=d(U^T)+O(d^2(U))$, in the third digit, too
little to be depicted.

GPUJACH1 timing started when all data were in place, and stopped when the
algorithm converged, but before the final collection of data.  The
similar holds for the sequential and MPI-parallel timings.  In all
cases, and for all algorithms, accumulation of $V^{-T}$ is included in
the timings.  If $V^{-T}$ is not needed (e.g., for the eigenvalue
problem), the speedups are a bit higher -- more than $4{\times}$ in
the 4-task case (for matrices of order $n \geq 4800$), and more than
$1.25{\times}$ in the 8-task task case (for matrices of order
$n \geq 8000$).  The effect of diagonal sorting on the speed of the
algorithm is also illustrated in Table~\ref{tbl:5.3}.  The similar
speedup occurs on the other problem sizes.

\begin{figure}[hbt]
\begin{center}
  \includegraphics{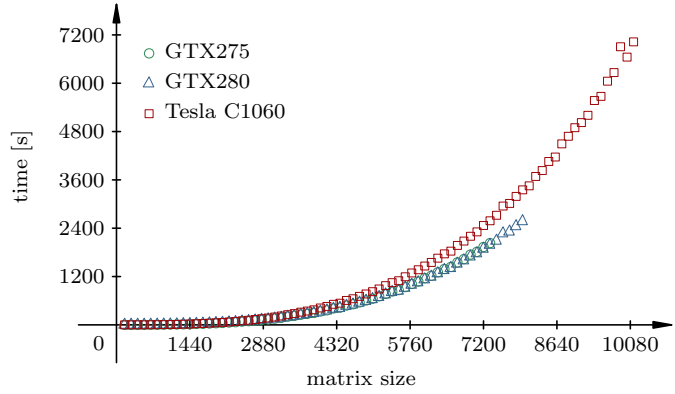}
\end{center}
\caption{Timing results of GPUJACH1 on square, full-rank matrices.}
\label{fig:5.2}
\end{figure}

Timing result on common test cases (Fig.~\ref{fig:5.2}) are perfectly
consistent between GTX275 and GTX280.  TeslaC1060 is a bit slower
because of the lower clock frequencies, but it is the only choice in
its generation for the large problems.

\begin{figure}[hbt]
\begin{center}
  \includegraphics{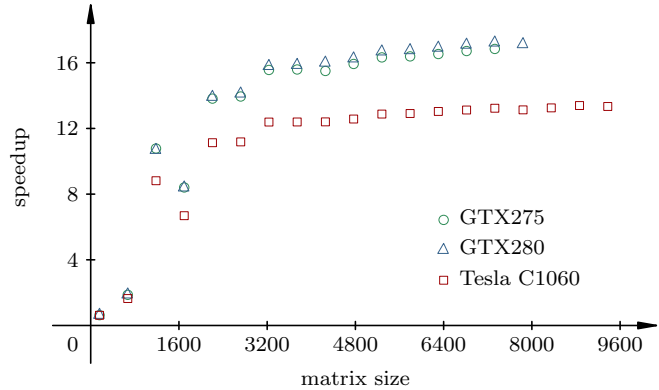}
\end{center}
\caption{Speedup of GPUJACH1 vs.~the sequential algorithm (with
  accumulation of $V^{-T}$).}
\label{fig:5.3}
\end{figure}

\begin{figure}[hbt]
\begin{center}
  \includegraphics{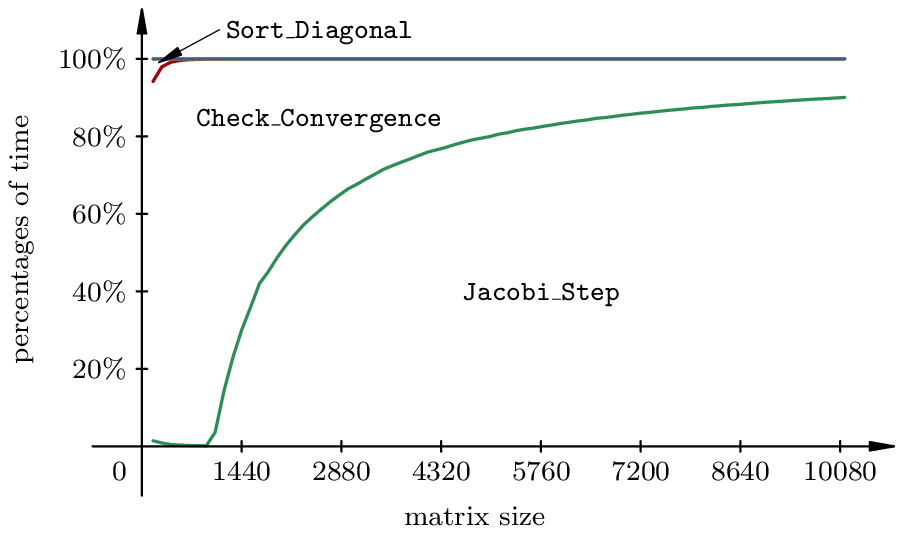}
\end{center}
\caption{Percentages of time.}
\label{fig:5.6}
\end{figure}

The speedup of GPUJACH1 vs.~the sequential algorithm is shown in
Fig.~\ref{fig:5.3}.  The inevitable overhead of convergence checking
(discussed in Section \ref{sec:4} and shown in Fig.~\ref{fig:5.6})
makes GPUJACH1 inappropriate for small matrices, but the speedup
quickly (at about $n=4800$) stabilizes up to $17{\times}$ in favor of
GPUJACH1.

\begin{figure}[hbt]
\begin{center}
  \includegraphics{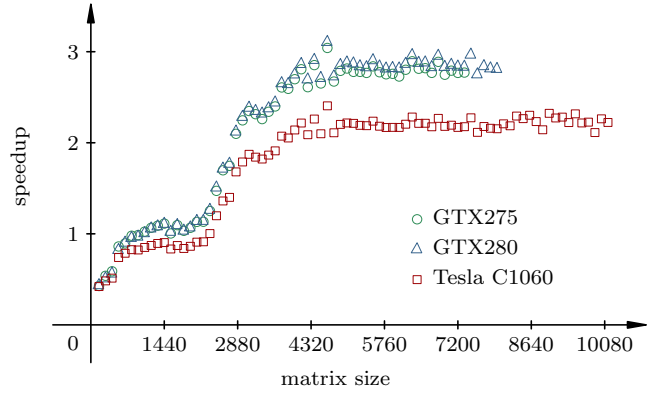}
\end{center}
\caption{Speedup of GPUJACH1 vs.~the 4-task parallel block-oriented
  algorithm (with accum.\ of $V^{-T}$).}
\label{fig:5.4}
\end{figure}

\begin{figure}[hbt]
\begin{center}
  \includegraphics{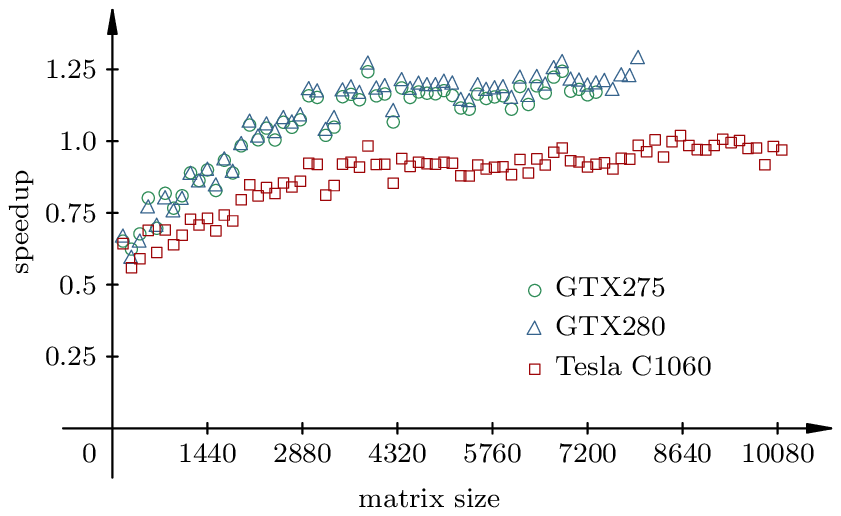}
\end{center}
\caption{Speedup of GPUJACH1 vs.~the 8-task parallel block-oriented
  algorithm (with accum.\ of $V^{-T}$).}
\label{fig:5.5}
\end{figure}

The speedup of GPUJACH1 vs.~the parallel block-oriented algorithm
(Figs.~\ref{fig:5.4} and \ref{fig:5.5}) is negligible (if none, i.e.,
CPU algorithm is faster), for matrices of order about $2000$.  One of
the reasons is efficient caching the blocking algorithms were designed
for, and GPUJACH1 had no caching opportunities whatsoever.  When the
CPU caches get too small for keeping the whole blocks without being
frequently evicted, GPUJACH1 attains it peak speedup.  The $3{\times}$
speedup over the 4-task case, compared to $17{\times}$ in the
sequential case (if accumulation of $V^{-T}$ is included) shows 
that the caching advantage of the block algorithm has dwindled.


GPUJACH1 algorithms pays off when \texttt{Jacobi\_Step} routine takes
a predominant part in the overall computation time.  The convergence
check is the main issue here, while the GPU sorting time is almost
negligible, except for the very small problems.  All the other
routines take well below 0.1\%\ of the time.  See Fig.~\ref{fig:5.6}
for details.

The profiling results of our algorithm are shown in
Table~\ref{tbl:5.4}.  The occupancy is not very high ($5$ blocks
residing on a SM, instead of theoretical maximum of $8$, because of the
register file starvation), which indicates that the future effort
should include simplifying the main \texttt{Jacobi\_Step} kernel,
maybe by breaking it down into a few lightweight kernels.  However,
the global memory throughput is satisfactory ($71.91\rm\%$ of the
theoretical limit of $102\rm\ GB/s$ on TeslaC1060).

\begin{table}[hbt]
  \begin{center}
    \begin{tabular}{@{}ccccccc@{}}
      \toprule
      computational     &  \multicolumn{3}{c}{global memory throughput [GB/s]} & instruction & achieved & register\\
      kernel            &  read & write & overall & throughput & occupancy & ratio\\
      \midrule
      \texttt{Jacobi\_Step}     & 50.0280 &     23.3173 & 73.3453 & 0.285392 & 0.3125 & 0.9375\\
      \texttt{Precompute\_Data} & 84.3262 & \hpz 0.3294 & 84.6556 & 0.219198 & 0.5000 & 0.5000\\
      \bottomrule
    \end{tabular}
  \end{center}
  \caption{Throughput and occupancy of HSVD on TeslaC1060 for matrices of order $4096$
    with $2048$ positive/negative signs in $J$ (with accumulation of $V^{-T}$).}
  \label{tbl:5.4}
\end{table}

We have also adapted our algorithm to be able to run more (but still
even) number of warps per block and to avoid the diagonal packing of
$(d,\rho,j)$ (see Fig.~\ref{fig:4.1} and discussion in Section
\ref{sec:4}).  However, the device code became more complex (i.e.,
uses $53$ registers and a couple of instructions more than GPUJACH1),
and the timing results were disappointing (see Table~\ref{tbl:5.5}).

\begin{table}[hbt]
  \begin{center}
    \begin{tabular}{@{}ccccc@{}}
      \toprule
      number of warps &   2 &   4 &   6 &   8 \\
      \midrule
      time            & 394 & 429 & 589 & 544 \\
      \bottomrule
    \end{tabular}
  \end{center}
  \caption{Speed of modified version of HSVD (without diagonal
    packing) on GTX275 for matrices of order $4096$ with $2048$
    positive signs in $J$ (with accum.\ of $V^{-T}$) and different
    number of warps per block.}
  \label{tbl:5.5}
\end{table}

A plausible explanation would be that the scheduling algorithm of
NVIDIA GPUs works reasonably well, and that it puts as much blocks on
a SM in a given time as possible.  Also, the blocks completes
execution when its slowest pair of warps does.  If you have a pair of
warps performing a trigonometric transformation, and another pair
performing a hyperbolic one in the same block, the hyperbolic one will
slow down the entire block.  As a part of a block cannot be replaced
on SM by a part of another block (only entire blocks are changeable),
it seems reasonable to have only one pair of warps in a block, which
turned out to be true.  Therefore, having more pairs of warps per
block could be viable if the computation of transformations, which
causes imbalance of execution time, is isolated into a separate,
lightweight kernel.

%
%%%%%%%%%%%%%%%%%%%%%%%%%%%%%%%%%%%%%%%%%%%%%%%%%%%%%%%%%%%%%%%%%%%%%%%%%%%%%%
%
\section{Conclusions and future work}
\label{sec:6}
%
%%%%%%%%%%%%%%%%%%%%%%%%%%%%%%%%%%%%%%%%%%%%%%%%%%%%%%%%%%%%%%%%%%%%%%%%%%%%%%
%
GPUJACH1 is fast and usable in its own right.  However, it can be incorporated
into the hybrid CPU--GPU parallel Jacobi framework with ease.  In
\cite{Singer-et-al-10}, besides the block-oriented, the parallel full-block
one-sided hyperbolic Jacobi algorithm is developed.  The full-block algorithm
diagonalizes its pivot block, while block-oriented only annihilates the
off-diagonal elements of a block once.  Instead of sequential diagonalization
of a block in each MPI process, GPUJACH1 can be plugged in as a direct
replacement.  The speedup should be the same as the speedup of GPUJACH1
vs.~the sequential algorithm, if each MPI process has access to its own GPU
unit and block sizes are large enough.

\looseness=-1
On GT200 chips there is no cache for the global memory, and the shared memory
is small.  Thus, it is difficult to avoid the slow BLAS 1 operations, and no
easy blocking is possible.  On the other hand, Fermi chips have multiple
levels of cache, and a larger shared memory.  That gives opportunity to employ
BLAS 3 operations for factorization in GPU and the shared-memory block
diagonalization (the full-block way).

Our future work also includes the symmetric indefinite factorization for GPU,
which is the key missing part to have a complete GPU-based symmetric
indefinite Jacobi-type eigensolver.

%
%%%%%%%%%%%%%%%%%%%%%%%%%%%%%%%%%%%%%%%%%%%%%%%%%%%%%%%%%%%%%%%%%%%%%%%%%%%%%%
%
\appendix
\section{Appendix}
\label{sec:7}
%
%%%%%%%%%%%%%%%%%%%%%%%%%%%%%%%%%%%%%%%%%%%%%%%%%%%%%%%%%%%%%%%%%%%%%%%%%%%%%%
%
There are many parallel strategies worth to study. A bunch of them are
described in \cite{Sameh-71,LukF-Park-87,LukF-Park-89}.  An obvious
choice of strategy is to take as many as possible, i.e.,
$\lfloor n / 2 \rfloor$ independent pivot pairs for annihilation in
each step, organized in $n - 1$ steps in a sweep for even $n$ and $n$
steps for odd $n$.  For many strategies which satisfy this property,
the proof of the convergence is not known.  Therefore, we loosen the
requirement on the maximal number of independent pivot pairs, but
still firmly require the convergence of the algorithm. 

The modulus strategy (in two different shifted forms) was described in
\cite{Sameh-71} and \cite{Hari-Veselic-87}. The name modulus strategy
for the first time appears in \cite{LukF-Park-89}, where on small
examples, its equivalence to the row cyclic strategy is illustrated.

Here we prove that our form of modulus strategy is weakly equivalent
to row cyclic strategy, and thus convergent.  The convergence of the
row cyclic strategy can be found in \cite[Theorem 2.3]{Veselic-93}.

Let us introduce a notation, taken from
\cite{Hari-SingerSanja-SingerSasa-10}.  In a pivot strategy, let $I(i, j)$
be the index at which the pair $(i, j)$ occurs.  Adjacent pivot pairs
$(i, j)$ and $(p, q)$ can swap their positions if $\{ i, j \} \cap \{
p, q \} = \emptyset$.  This transposition of pairs is called
\emph{admissible transposition\/}. Two pivot strategies $O$ and $O'$
are \emph{equivalent\/} if one can be transformed to the other by a
finite number of admissible transpositions.

Two pivot strategies $O$ and $O'$ are \emph{shift equivalent\/}, if
one is obtained form the other by cyclic shift of the pivot pairs,
i.e., if the position of each pair in strategy $O$ is $I(i, j)$, its
new position in strategy $O'$ is
\begin{displaymath}
  I'(i, j) := (I(i, j) + c) \bmod n_p,
\end{displaymath}
where $n_p = n (n - 1) / 2$ and $c \in \mathbb{Z}$.

Two pivot strategies $O$ and $O'$ are \emph{weakly equivalent\/}, if
there exist strategies $O_i$, $1 \leq i \leq r$, for some $r \geq 1$,
such that in the sequence of strategies $O, O_1, \ldots, O_r, O'$ each
two adjacent terms are either equivalent or shift equivalent.
Obviously, proof of convergence for one strategy in the sequence is
sufficient to prove that all strategies in the sequence are
convergent.

The convergence of an algorithm is usually proved for the two-sided
algorithm.  Since the one-sided transformations on a factor of a
matrix are in theory (need not be numerically) the same as the
two-sided algorithm, we immediately have the proof of convergence for
the one-sided algorithm.

First we prove that the antidiagonal strategy is equivalent to
row cyclic strategy, and second that the modulus strategy is weakly
equivalent to antidiagonal strategy.

\begin{figure}[hbt]
\begin{center}
  \includegraphics{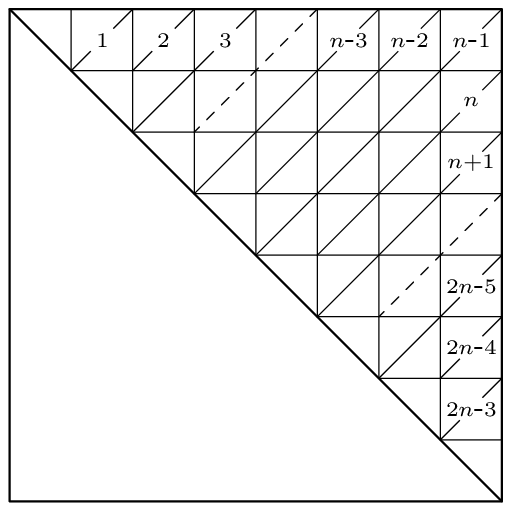}\qquad
  \includegraphics{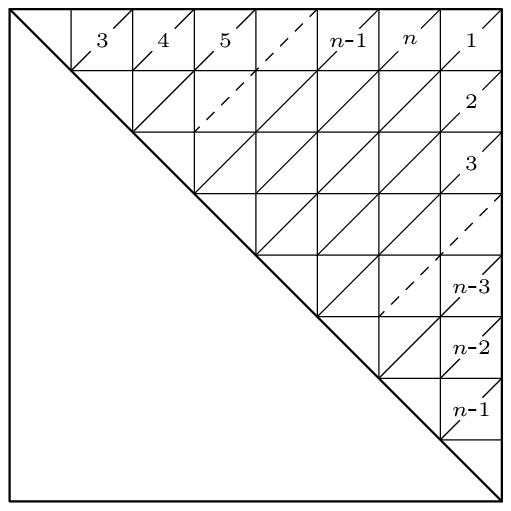}
\end{center}
\caption{Antidiagonal strategy (left), and the modulus strategy
  (right).  The order of steps are labeled on each antidiagonal.}
\label{fig:A.1}
\end{figure}

Antidiagonal strategy consists of the sequence of steps shown in
Fig.~\ref{fig:A.1} (left).  The steps of the modulus strategy are
shown in Fig.~\ref{fig:A.1} (right).

Table \ref{tbl:A.1} shows the order of annihilation in the antidiagonal
strategy.
\begin{table}[hbt]
  \begin{center}
      \begin{tabular}{@{}cccccccc@{}}
        \toprule
        step     & \multicolumn{7}{c}{simultaneously annihilated pivot pairs} \\
        \midrule
        1        & $(1, 2)$ &            &            &          &            &  & \\
        2        & $(1, 3)$ &            &            &          &            &  & \\
        3        & $(1, 4)$ & $(2, 3)$   &            &          &            &  & \\
        4        & $(1, 5)$ & $(2, 4)$   &            &          &            &  & \\
        5        & $(1, 6)$ & $(2, 5)$   & $(3, 4)$   &          &            &  & \\
        $\vdots$ & $\vdots$ &            & $\vdots$   &          &            &  & \\
        $n - 1$  & $(1, n)$ & $(2, n-1)$ & $(3, n-2)$ & $\cdots$ & $(k, \ell)$ & & \\
        $n$      &          & $(2, n)$   & $(3, n-1)$ & $\cdots$ &
        $\begin{cases}
          (k + 1, \ell), & \text{$n$ odd} \\
          (k, \ell + 1), & \text{$n$ even}
        \end{cases}$ & & \\
        $\vdots$ &          &            & $\ddots$   &          & $\ddots$ &           &  \\
        $2n - 2$ &          &            &            &          & & $(n-2, n-1)$ & $(n-2, n)$ \\
        $2n - 3$ &          &            &            &          & &           & $(n-1, n)$ \\
        \bottomrule
      \end{tabular}
  \end{center}
\caption{Order of annihilation in the antidiagonal strategy.}
\label{tbl:A.1}
\end{table}

There is a minor difference between odd and even $n$ in the step $n-1$
in Table \ref{tbl:A.1}, i.e.,
\begin{displaymath}
  k = \begin{cases}
    \frac{n-1}{2}, & \text{$n$ odd}, \\
    \frac{n}{2},   & \text{$n$ even},
  \end{cases}
  \qquad
  \ell = \begin{cases}
    k + 2, & \text{$n$ odd}, \\
    k + 1, & \text{$n$ even}.
  \end{cases}
\end{displaymath}

By admissible transpositions, pivot pairs from the first column
(pairs $(1, j)$) in Table \ref{tbl:A.1} could be written
before all pairs $(2, {} \cdot {}), (3, {} \cdot {}),
\ldots, (k - 1, {} \cdot {})$, since pivot indices $(1, j)$, $j = 5,
\ldots, n$ are disjoint with indices of these columns (the first index
in $(1, j)$ is always smaller than the first index in other columns,
and the second is always greater than the second index in other
columns).  The rest of the proof, for indices in the second, third,
and other columns follows by induction over column index.  This
completes the proof of equivalence of the antidiagonal and the
row cyclic strategy.

The second step in the proof is to show that antidiagonal and the
modulus strategy are weakly equivalent. 

According to Fig.~\ref{fig:A.1} modulus strategy consists of
either one or two steps of the antidiagonal strategy.  Note that the
antidiagonal strategy, which consists of steps
\begin{displaymath}
  O = (1, 2, 3, \ldots, n-2, n-1, n, n+1, n+2, \ldots, 2n-4, 2n-3),
\end{displaymath}
is shift equivalent to strategy which consists of steps
\begin{displaymath}
  O_1 = (n-1, n, n+1, n+2, \ldots, 2n-4, 2n-3, 1, 2, 3, \ldots, n-2).
\end{displaymath}
The pivot indices in step $1$ (pair $(1, 2)$) are disjoint with 
the indices contained in steps $n+2, \ldots, 2n-3$ since the first index
is at least $3$, and the second is at least $\ell+1 > 2$.  Thus,
cyclic strategy $O_1$ is equivalent to
\begin{displaymath}
  O_2 = (n-1, n, n+1, 1, n+2, \ldots, 2n-4, 2n-3, 2, 3, \ldots, n-2).
\end{displaymath}
The similar reasoning holds for step $2$, with indices disjoint with
indices in steps $n+3, \ldots, 2n-3$, indices from step $3$ are
disjoint with indices in steps $n+4, \ldots, 2n-3$, and so on.

The final sequence of steps is
\begin{displaymath}
  O' = (n-1, n, n+1, 1, n+2, 2, \ldots, 2n-3, n - 2),
\end{displaymath}
which is in fact modulus strategy, since step $n-1$ of the antidiagonal
strategy is the first step of the modulus strategy, step $n$ is the
second step of the modulus strategy, the third step of the modulus strategy
consists of the steps $n+1, 1$ from the antidiagonal strategy, and so on. 

This proves that the modulus strategy is weakly equivalent to
row cyclic strategy.

%
%%%%%%%%%%%%%%%%%%%%%%%%%%%%%%%%%%%%%%%%%%%%%%%%%%%%%%%%%%%%%%%%%%%%%%%%%%%%%%
%
\subsection*{Acknowledgements}
%
%%%%%%%%%%%%%%%%%%%%%%%%%%%%%%%%%%%%%%%%%%%%%%%%%%%%%%%%%%%%%%%%%%%%%%%%%%%%%%
%
We would like to thank Prof.~Većeslav Čorić from Faculty of Mechanical
Engineering and Naval Architecture, University of Zagreb, for helping us
acquiring a part of the testing equipment, and Prof.~Vjeran Hari from
the Department of Mathematics, University of Zagreb for helping us
with the proof of the convergence of the modulus pivot strategy.
We also thank the anonymous referees for their helpful suggestions, which
substantially improved the paper.
%
%%%%%%%%%%%%%%%%%%%%%%%%%%%%%%%%%%%%%%%%%%%%%%%%%%%%%%%%%%%%%%%%%%%%%%%%%%%%%%
%
%\bibliographystyle{spmpsci}
%\bibliography{ref}

\end{document}